\documentclass[11pt]{article}

\usepackage{graphicx}
\usepackage[top=30mm, bottom=30mm, left=40mm, right=40mm]{geometry}
\usepackage{amsmath}
\usepackage{mathtools}
\DeclareMathOperator{\wt}{wt}

\begin{document}


\title{Ulam-Warburton Automaton - Counting Cells with Quadratics}  
\author{Mike Warburton}
\date{\today}
\maketitle

\begin{abstract}
This paper is about a sequence of quadratic functions that enumerate the total number of \textsc{on} cells up to and including generation $n$ of the Ulam-Warburton cellular automaton, where $n$ has the form $n_m=m\cdot2^k$.
\end{abstract}

		Keywords: cellular automata (CA), enumeration, 		Ulam-Warburton, UWCA.

\section{Introduction}

The origins of cellular automata go back to Stanislaw Ulam in 1929 \cite{UW1}, he later explored these ideas with J. C. Holladay and Robert Schrandt \cite{UW2}. The first function counting the total number of \textsc{on} cells is the quadratic for the sharp upper bound occurring at generations $n = 2^k$ \cite{UW3},  this sparked interest and the fractal like object was named the Ulam-Warburton Cellular Automaton (UWCA) in 2003 \cite{UW4}. Since then mathematicians have connected the UWCA with various objects including the Toothpick sequence \cite{UW5}, Nim Fractals \cite{UW6} and the Sierpinski triangle \cite{UW7}.

\begin{center}
\includegraphics[scale=.7]{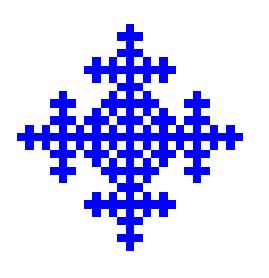} 
\end{center}

Figure 1: The Ulam-Warburton cellular automaton  ($n=14$)\\

The Ulam-Warburton cellular automaton is a 2-dimensional fractal pattern that grows on a grid of cells consisting of squares. Starting with one square initially \textsc{on} and all others \textsc{off}, successive iterations are generated by turning \textsc{on} all squares
that share precisely one edge with an \textsc{on} square. This is the von Neumann neighbourhood.

\noindent Starting the generation index when the first cell is \textsc{on}, the sharp upper bound function  agrees with  the total number of \textsc{on} cells when $n=2^k$ and is described by.

\begin{equation*} 
 U_{sub}(n)=\frac{4}{3} n^2-\frac{1}{3}.
\end{equation*}

\noindent This paper extends this result to all generations of the UWCA.

\section{The Development of the Quadratics}

Let $u(n)$ denote the number of \textsc{on} cells at the $n^{th}$ stage.\\ 
$u(0)=0,u(1)=1$ and for $n \geq 2$.

\begin{equation*} 
 u(n)= \frac{4}{3}3^{\wt(n-1)}.
\end{equation*}

\noindent Where $\wt(n)$ is the Hamming weight function which counts the number of 1's in the binary expansion of $n$ \cite{UW5}.\\

\noindent Let $U(n)$ denote the total number of \textsc{on} cells after $n$ stages.

\begin{equation*} 
 	U(n)=\sum_{i=0}^n u(i),
\end{equation*}

\begin{equation*} 
 	U(n)=\frac{4}{3}\sum_{i=0}^{n-1}3^{\wt(i)}-\frac{1}{3}.
\end{equation*}

\noindent We now consider integer sequences $n_m$ based on the form $n_m=m\cdot2^k$ where $m \geq 1$ and $k \geq 0$. The total number of \textsc{on} cells for these sequences becomes.

\begin{equation}  
 	U_m(n_m)=\frac{4}{3}\sum_{i=0}^{m2^k-1}3^{\wt(i)}-\frac{1}{3}.
\end{equation}

\noindent Using this notation the expression for the sharp upper bound $U_{sub}(n)$ becomes  $U_1({n_1})$ and therefore.

\begin{equation} 
 	\frac{4}{3}\sum_{i=0}^{2^k-1}3^{\wt(i)}-\frac{1}{3}=\frac{4}{3} n_1^2-\frac{1}{3}.
\end{equation}

\noindent Which is the first of a sequence of quadratics, and in terms of $k$ we have.

\begin{equation*}  
	U_1(k)=\frac{4}{3} 2^{2k}-\frac{1}{3}.
\end{equation*}

\noindent Returning to equation (1) for $U_m(n_m)$ and introducing the relationship.

\begin{equation} 
 	\frac{4}{3}\sum_{i=0}^{m2^k-1}3^{\wt(i)}-\frac{1}{3}=	\sum_{i=0}^{m-1}3^{\wt(i)}\frac{4}{3}\sum_{i=0}^{2^k-1}3^{\wt(i)}-\frac{1}{3}.
\end{equation}

\noindent Let

\begin{equation*} 
	a_m= \sum_{i=0}^{m-1}3^{\wt(i)}
\end{equation*}

\noindent This is OEIS \cite{UW8} sequence A130665. \\

\noindent Substituting $a_m$ in to equation (3) and using equation (1) we have.

\begin{equation*} 
 	U_m(n_m)=a_m\frac{4}{3}\sum_{i=0}^{2^k-1}3^{\wt(i)}-\frac{1}{3}.
\end{equation*}

\noindent As $n_m=m\cdot2^k = m\cdot n_1$ therefore $n_1^2=\frac{n_m^2}{m^2}$ using this with equation (2) we have.

\begin{equation} 
 	U_m(n_m)=\frac{a_m}{m^2} \frac{4}{3}n_m^2 - \frac{1}{3}.
\end{equation}

\noindent This is the result we were aiming for and in terms of k we have.

\begin{equation*} 
 	U_m(k)=a_m \frac{4}{3}2^{2k} - \frac{1}{3}.
\end{equation*} \\

\begin{tabular}{|r|r|r|r|r|r|r|r|r|}
\hline
$k$& $n_1$& $U_1$& $n_3$ &$U_3$& $n_5$& $U_5$& $n_7$& $U_7$ \\ \hline
0 &1  &1  &3   &9  &5  &25  &7  &49 \\ \hline
1 &2  &5  &6  &37 &10 &101 &14 &197 \\ \hline
2 &4 &21 &12 &149 &20 &405 &28 &789 \\ \hline
3 &8 &85 &24 &597 &40 &1,621 &56 &3,157 \\ \hline
4 &16 &341 &48 &2,389 &80 &6,485 &112 &12,629 \\ \hline
5 &32 &1,365 &96 &9,557 &160 &25,941 &224 &50,517 \\ \hline
6 &64 &5,461 &192 &38,229 &320 &103,765 &448 &202,069 \\ \hline
7 &128 &21,845 &384 &152,917 &640 &415,061 &896 &808,277 \\ \hline
8 &256 &87,381 &768 &611,669 &1,280 &1,660,245 &1792 &3,233,109 \\ \hline

\end{tabular}\\
\\  

\noindent Table 1: The odd numbered sequences $U_1$ to $U_7$ that give the total number of \textsc{on} cells in the sequences $n_1$ to $n_7$.

\section{The Limit Inferior and Limit Superior}

We have \cite{UW9}.

\begin{equation*}  
 	0.9026116569...= \liminf\limits_{n\to \infty} \frac{U(n)}{n^2} < \limsup\limits_{n\to \infty} \frac{U(n)}{n^2} =\frac{4}{3}.
\end{equation*} \\

\noindent We can see from equation (4) that the coefficient of the quadratic term is the dominant expression in $\frac{U(n)}{n^2	}$ and is a rational number.\\
\\


\end{document}